\documentclass[12pt]{amsart} 

\title{The stable mapping class group of simply connected 4-manifolds} 

\author{Jeffrey Giansiracusa} 
\address{I.H.E.S. \\ Le Bois-Marie \\ 35 Rue de Chartres\\
Bures-sur-Yvette F-91440 \\ France} 
\email{giansira@maths.ox.ac.uk}

\usepackage[margin=3.0cm,dvips]{geometry}

\usepackage{pst-all}
\usepackage{mathrsfs}
\usepackage{amssymb}
\usepackage[cmtip,arrow]{xy}
\usepackage{pb-diagram,pb-xy}

\newtheorem{theorem}{Theorem}[section] 
\newtheorem{lemma}[theorem]{Lemma} 
\newtheorem{proposition}[theorem]{Proposition}
\newtheorem{corollary}[theorem]{Corollary}

\theoremstyle{remark} 
\newtheorem{remark}[theorem]{Remark}
\newtheorem{definition}[theorem]{Definition}

\newcommand{\Diff}{\mathrm{Diff}}
\newcommand{\Emb}{\mathrm{Emb}}
\newcommand{\Aut}{\mathrm{Aut}}

\DeclareMathOperator*{\colim}{colim}
\newcommand{\C}{\mathbb{C}}

\newcommand{\N}{\mathbb{N}}
\newcommand{\Z}{\mathbb{Z}}
\newcommand{\Q}{\mathbb{Q}}

\newcommand{\fix}{\mbox{ fix }}
\newcommand{\calC}{\mathcal{C}}
\newcommand{\calK}{\mathcal{K}}
\newcommand{\CPCPB}{\mathbb{C} P^2 \# \overline{\mathbb{C} P^2}}

\begin{document}
\vspace*{-1cm}
\begin{abstract}
  We consider mapping class groups $\Gamma(M) = \pi_0 \Diff(M \fix
  \partial M)$ of smooth compact simply connected oriented
  4--manifolds $M$ bounded by a collection of 3--spheres.  We show
  that if $M$ contains $\C P^2$ or $\overline{\C P}^2$ as a connected
  summand then all Dehn twists around 3--spheres are trivial, and
  furthermore, $\Gamma(M)$ is independent of the number of boundary
  components.  By repackaging classical results in surgery and
  handlebody theory from Wall, Kreck and Quinn, we show that the
  natural homomorphism from the mapping class group to the group of
  automorphisms of the intersection form becomes an isomorphism after
  stabilization with respect to connected sum with $\C P^2 \#
  \overline{\C P^2}$.  We next consider the 3+1 dimensional cobordism
  2--category $\calC$ of 3--spheres, 4--manifolds (as above) and
  enriched with isotopy classes of diffeomorphisms as 2--morphisms.
  We identify the homotopy type of the classifying space of this
  category as the Hermitian algebraic $K$-theory of the integers.  We
  also comment on versions of these results for simply connected spin
  4--manifolds.  Finally, we observe that a related 4--manifold operad
  detects infinite loop spaces.
\end{abstract}
\maketitle

\section{Introduction}

In this paper we shall be interested in the mapping class groups
\[
\Gamma(M) := \pi_0 \Diff(M \fix \partial M)
\]
of smooth compact simply connected oriented 4-manifolds $M$ bounded by
any number of ordinary 3--spheres.  Our strategy is to compare the
mapping class group to the group $\Aut(Q_M)$ of automorphisms of the
intersection form $Q_M$ of $M$, which is an object with a more
algebraic character and which has a far clearer structure.  This
problem is analogous to the well-studied problem in surface theory of
understanding the Torelli group.  For a surface of genus $g$, the
Torelli group $I_g$ is the kernel of the natural map $\Gamma_g \to
Sp_{2g}\Z$ which sends an isotopy class to its induced automorphism of
the intersection form (which is a symplectic form for surfaces).  The
Torelli group for surfaces contains quite a lot of rich structure,
even in the stable setting of infinite genus.  In contrast we show
that the situation is as different as possible for simply connected
4-manifolds.  After stabilizing by taking connected sums with $\C P^2
\# \overline{\C P^2}$, the stable mapping class group becomes
isomorphic to the stable automorphism group, and thus the stable
Torelli group of simply connected 4-manifolds vanishes.

We next turn towards cobordism categories of simply connected
4-manifolds.  This is motivated by Witten's and Morava's ideas about
``topological quantum gravity'' (e.g. \cite{Witten},
\cite{Morava-grav}, \cite{Morava-tate}), and also by results in
dimension 2 (such as the proof of the stable Mumford conjecture) which
have origins tracing back to Tillmann's analysis of the analogous
category in dimension 2 \cite{Tillmann-mcg-homotopy},
\cite{Tillmann-splitting}.  We find that the homotopy type of the
classifying space of the cobordism 2--category of \emph{(3--spheres,
  simply connected 4-manifolds, and isotopy classes of
  diffeomorphisms)} is equivalent as infinite loop spaces to the
Hermitian algebraic $K$-theory of the integers.  We also give a
variant of this analysis for the restriction to spin 4-manifolds.

Finally, we observe that the operad constructed from mapping
class groups of simply connected 4--manifolds fits into Tillmann's
framework \cite{Tillmann-operad}.  Thus the 4--manifold operad also
detects infinite loop spaces.

\subsection*{Acknowledgements}

I would like to thank Ralph Cohen for communicating to me the
unpublished results obtained by his student, Jun Yan, which motivated
this project and overlap with some of the results obtained in this
paper.  I would also like to thank my advisor, Ulrike Tillmann, for
countless conversations, and John Rognes and an anonymous referee for
carefully reading the manuscript and pointing out several corrections.
I gratefully acknowledge the support of an NSF Graduate Research
Fellowship.

\section{Overview and statement of results}

Given a simply connected 4-manifold $M$, we let 
\[ \Gamma(M) := \pi_0 \Diff(M \fix \partial M) \] denote the mapping
class group of isotopy classes of orientation-preserving
diffeomorphisms which are the identity on $\partial M$.  The
intersection pairing $Q_M: H_2 (M;\Z)\otimes H_2 (M;\Z) \to \Z$
(thought of as a \emph{symmetric bilinear form}) determines a group
$\Aut(Q_M) \subset GL(H_2 (M;\Z))$ of automorphisms which preserve the
intersection pairing.  There is a homomorphism $\Gamma(M) \to
\Aut(Q_M)$ induced by sending a diffeomorphism to the induced
automorphism on $H_2$, and the kernel of this map is defined to be the
4-manifold Torelli group $I(M)$.

Let us briefly recall some known facts about this map for surfaces and
4-manifolds.  For surfaces, as soon as the genus is larger than 1 the
components of the diffeomorphism group are contractible so $\Diff(F_g)
\to \Gamma_g$ is a homotopy equivalence.  The rational cohomology of
$\Gamma_g$ (in a stable range proportional to the genus) is a
polynomial algebra $\Q[\kappa_1, \kappa_2, \ldots]$ on the
Miller-Morita-Mumford classes $\kappa_i$.  The odd half of these
classes pull-back from the integral symplectic group via the map
$\Gamma_g \to Sp_{2g}(\mathbb{Z})$, but the even half do not. Even at
the cohomological level the two groups are fairly different.

However, the situation is markedly simpler for \emph{topological}
4-manifolds.
\begin{theorem}[\cite{Quinn}]
For a simply connected topological 4-manifold $M$,
\[\pi_0 \mathrm{Homeo}(M) \cong \Aut(Q_M).\]
\end{theorem}
Of course, once we set foot in the land of smooth manifolds, the
situation becomes somewhat more interesting, for Ruberman
\cite{Ruberman} has constructed examples of smooth 4-manifolds for
which the homomorphism $\pi_0 \Diff(M) = \Gamma(M) \to \Aut(Q_M)$ has
non-finitely generated kernel!  One detects and distinguishes elements
in the kernel using a gauge theoretic invariant.  Note that in a rough
sense gauge theory tends to only detect properties which are
\emph{unstable} with respect to connected sum.  For instance, the
Donaldson polynomial vanishes after a single connected sum with $S^2
\times S^2$ \cite[Theorem 1.3.4, p.26]{Donaldson-Kronheimer}.
Furthermore, if we allow ourselves to start taking connected sums with
$S^2\times S^2$ then Wall's stable classification theorem \cite{Wall2}
tells us that the stable diffeomorphism type is determined entirely by
the intersection form.  Motivated by these examples of how the unusual
phenomena of smoothness in dimension 4 tend to go away in
stabilization, there is hope that we may understand $\Gamma(M)$
\emph{stably}.

In the world of surfaces stabilization is a familiar concept---one
stabilizes by sequentially gluing on tori to let the genus go to
infinity.  There are a multitude of results which illustrate the
utility of considering this stabilization.  Tillmann's theorem
\cite{Tillmann-mcg-homotopy} that the classifying space of the stable
mapping class group is an infinite loop space after applying Quillen's
plus construction is one such example.  This and many other results
stem from Harer's homological stability theorem \cite{Harer}, which is
considered to be one of the high-points of surface theory.  It says
that increasing the genus and number of boundary components induces
isomorphisms on the homology of the mapping class groups in a stable
range of degrees that depends only on the genus.  (The stability range
was later improved by Ivanov \cite{Ivanov}.)

For 3--manifolds, Hatcher and Wahl \cite{Hatcher-Wahl} have proven that
the homology of the mapping class group modulo all Dehn twists around
embedded 2--spheres has a similar homological stability property.
(Note that in dimension 2 Dehn twists are known by the Dehn-Lickorish
Theorem \cite{Dehn, Lickorish} to generate the mapping class group.)
In dimension 4 however there is no known theorem analogous to Harer's
stability.  Nevertheless, we shall find that there is still
appreciable utility in studying 4-manifold mapping class groups under
stabilization.

The particular stabilization we focus on in the present paper is that
of repeatedly taking connected sums with $\mathbb{C}P^2 \#
\overline{\mathbb{C}P^2}$.  (Here $\overline{\C P^2}$ denotes the
complex projective plane with the opposite orientation.)  Though
perhaps less familiar than using $S^2\times S^2$, this stabilization
makes sense for the following two reasons. (i) There is a
diffeomorphism $(S^2\times S^2) \# \overline{\mathbb{C}P^2} \cong
\mathbb{C}P^2 \# \overline{\mathbb{C}P^2} \#
\overline{\mathbb{C}P^2}$, so our stabilization process automatically
implicitly contains the more familiar stabilization with respect to
$S^2\times S^2$.  (ii) Intersection forms (which are integral
unimodular forms) are either \emph{even} or \emph{odd}, and
\emph{definite} or \emph{indefinite}.  Since connected sum of
manifolds corresponds to direct sum of intersection forms, the odd
indefinite quadrant of the classification is the only one which cannot
be exited by taking connected sums, and a form in this quadrant is
isomorphic to $n(1)\oplus m(-1)$ by the classical Hasse-Minkowski
classification.  The intersection form of $\mathbb{C}P^2 \#
\overline{\mathbb{C}P^2}$ is $(1)\oplus (-1)$, so our stabilization
process puts us immediately into the land of odd indefinite forms and
sends the numbers of $(1)$ and $(-1)$ summands both to infinity.

To stabilize the mapping class group one must have a way to extend an
isotopy class $[\phi]$ across a connected sum.  In general this is
impossible since one needs a fixed disc in which to perform the
cutting and pasting, and a chosen representative $\phi$ need not fix a
disc anywhere.  Of course, one can always choose a representative
which does fix a disc, but then the isotopy class resulting from
extending this representative may depend on the choice.  Instead we
shall use manifolds with boundary and stabilize by gluing along the
boundary.  Let $M$ be a 4-manifold bounded by some number of ordinary
3--spheres, and let $X$ denote $\mathbb{C}P^2
\#\overline{\mathbb{C}P^2}$ with the interiors of two discs removed.
We may glue $X$ along a selected boundary component of $M$ to obtain a
new manifold denoted $MX$, and then iterate by gluing along the
remaining boundary of $X$, as in Figure \ref{stab-fig}.  Extension by
the identity on $X$ determines a system of maps
\[\Gamma(M) \to \Gamma(MX) \to \Gamma(MX^2) \to \cdots \]
and the \emph{stable mapping class group of $M$}, written
$\Gamma_\infty(M)$, is defined to be the colimit of this sequence.
\begin{figure}
  \caption{Stabilizing by gluing on copies of $X = \CPCPB$. \label{stab-fig}}
  \begin{center}
    \vspace{5mm}
    \input{stab-fig.pst}
  \end{center}
\end{figure}

On the intersection form, gluing $M$ to $N$ along a boundary 3--sphere
(or even just a homology 3--sphere) induces an inclusion $\Aut(Q_M)
\hookrightarrow \Aut(Q_{MN})$ coming from the block addition of
intersection forms, $Q_M \mapsto Q_{MN} = Q_M
\oplus Q_N$.  We thus define the \emph{stable automorphism group of
  $M$},
\[\mathrm{Aut}_\infty(M) := \colim_{n\to\infty} \{ \Aut(Q_M) \hookrightarrow
\Aut(Q_{MX}) \hookrightarrow \cdots \}.\]

Our main result is the following.
\begin{theorem}\label{stable-isomorphism}
  The stable groups $\Gamma_\infty(M)$ and $\Aut_\infty(Q_M)$ do not
  depend on the choice of the initial manifold $M$ within the class of
  smooth compact oriented simply connected 4-manifolds bounded by a
  collection of ordinary 3--spheres, and in particular, they are
  independent of the number of boundary components of $M$.
  Furthermore, $\Gamma_\infty \cong \Aut_\infty \cong
  O_{\infty,\infty}(\Z)$, with the first isomorphism induced by the
  natural map $\Gamma \to \Aut$ and the second coming from choosing a
  basis.
\end{theorem}
(Saeki \cite{Saeki} has independently proved a more general version of
this theorem by similar methods.)  Here $O_{\infty, \infty}(\Z)$ is
the group of automorphisms of the quadratic form $\infty(1)\oplus
\infty (-1)$ on $\Z^\infty$.  This group is closely related to the
Hermitian $K$-theory of the integers.  See the discussion immediately
after Theorem \ref{cat-equiv}.

\begin{corollary}
The 4-manifolds stable Torelli group $I_\infty$ is zero.
\end{corollary}

Along the way we find it necessary to analyze the elements of the
mapping class group represented by Dehn twists around 3--spheres
embedded with trivial normal bundle.  In section \ref{twist-section}
we define these elements and show that the set of twists around
boundary 3--spheres generates the kernel of the surjective
homomorphism induced by filling the boundary components with discs.
Furthermore, we find that a $\C P^2$ summand will kill all Dehn
twists.

\begin{theorem}\label{twists-are-trivial1}
  Let $M$ be a smooth compact 4-manifold (nontrivial fundamental group
  and nonempty boundary are allowed) of the form $N\# \C P^2$ (or $N\#
  \overline{\C P^2}$).  Then any Dehn twist on $M$ is isotopic to the
  identity.
\end{theorem}

\begin{corollary}\label{twists-are-trivial2}
  Let $M$ be a simply connected oriented \emph{closed} 4-manifold and
  let $M'$ be the result of removing the interiors of $n$ disjoint
  discs in $M$; let $K$ denote the kernel of the surjective map
  $\Gamma(M') \to \Gamma(M)$.  
\begin{enumerate}
\item[$(i)$] $K$ is generated by Dehn twists around the boundary spheres of $M'$.
\item[$(ii)$] If $M$ is of the form $N\#\C P^2$ (or $N\# \overline{\C P^2}$)
  then $K=0$.
\item[$(iii)$] If $M$ is spin then $K$ is either $(\Z/2)^{n-1}$ or $(\Z/2)^n$.
\end{enumerate}
\end{corollary}

The above corollary may be viewed as a very strong form of stability
with respect to increasing the number of boundary components for
mapping class groups of 4-manifolds.  It holds at the level of groups
and it requires only a single stabilization step, whereas for surfaces
Harer-Ivanov stability says that the analogue of the above map is
merely a homology isomorphism and only in a stable range depending on
the genus.  The proofs of \ref{twists-are-trivial1} and
\ref{twists-are-trivial2} are entirely elementary, whereas the proof
Harer-Ivanov stability requires the machinery designed to deal with
curve-complexes.

Theorem \ref{stable-isomorphism} is partly a repackaging of classical
theorems in 4-manifold topology due to Wall, Kreck, and Quinn.  Wall
\cite{Wall1} proved that $\Gamma(M)\to \Aut(Q_M)$ is surjective when
$M$ is indefinite and contains $S^2\times S^2$ as a connected summand.
On the other hand, Kreck \cite{Kreck} proved that the map is
\emph{always} injective, once one descends from isotopy to
pseudo-isotopy.  Finally, Quinn \cite{Quinn} proved that
pseudo-isotopy implies isotopy in a stable sense.  Together these
results yield a lifting of the automorphism group of the intersection
form of $M$ into the stabilized mapping class group of $M$ and in the
colimit this becomes the inverse to $\Gamma_\infty \to \Aut_\infty$.
This material is covered in more detail in section
\ref{stable-groups}.

Our results above have a close connection to what Morava
\cite{Morava-grav} calls a \emph{theory of topological gravity} in
4-dimensions.  Generalizing Witten's theory \cite{Witten} in
2--dimensions, Morava defines such a theory to be a representation of
a topological 2--category $\mathcal{G}$ where objects are
3--manifolds, morphisms are 4--dimensional cobordisms and 2--morphisms
are diffeomorphisms of the cobordisms. There are many possible
variations in the definition of this category.  The symmetric monoidal
product given by disjoint union implies that the classifying space of
the cobordism category is an infinite loop space, and a representation
induces an infinite loop map into some variant of a $K$--theory
infinite loop space.  Thus a theory of topological gravity produces an
element in some version of the $K$--theory of $B\mathcal{G}$.  As a
rough first step towards constructing or classifying topological
gravity theories one would thus like to understand the homotopy type
of $B\mathcal{G}$ (or at least some version of its $K$--theory).

The recent work of Galatius, Madsen, Tillmann and Weiss \cite{GMTW}
determines the homotopy types of many versions of the category
$\mathcal{G}$ in terms of more accessible spaces: the zeroth spaces of
certain Thom spectra.  Given their result, one might wonder why it
should be necessary to study the homotopy type with any other methods.
However, their argument only applies when $\mathcal{G}$ is maximal in
the sense that it is built using \emph{all} manifolds of appropriate
dimensions; their theorem does not determine the homotopy type of
sub-categories which are obtained by restricting the objects or
morphisms.

Morava \cite{Morava-tate} has indicated that one such subcategory out
of reach to the GMTW theorem is potentially interesting; he argues
based on an analogy with the Virasoro algebra that symmetries of the
Tate cohomology $t^*_{SU(2)}kO$ should play a role in the
representation theory of the cobordism category of smooth spin
4--manifolds bounded by ordinary 3--spheres.

In sections \ref{cats} and \ref{cats-proof} of this paper we study the
homotopy type of two simplified variants of Morava's category.  Let
$\calC$ denote the cobordism 2--category where the objects of $\calC$
are unions of 3--spheres, the morphisms are disjoint unions of simply
connected ``tree-like'' 4--manifolds --- meaning that each component
has precisely one outgoing boundary sphere; this is imposed so that
compositions stay within the simply connected realm --- and the
2--morphisms are \emph{isotopy classes} of diffeomorphisms.  

We will construct a map from $\Omega B\calC$ into the Hermitian algebraic
$K$--theory space $\Z^2 \times BO_{\infty,\infty}(\Z)^+$; it is induced
by the natural 2--functor from $\calC$ to the 2--category $\calK$
constructed similarly, with 2--morphisms now being isometries (with
respect to the intersection forms) of $H_2$.  We define these
categories more carefully in section \ref{cats}.  These categories are
strict symmetric monoidal under disjoint union, and hence their
classifying spaces are infinite loop spaces.  Our main result here is:

\begin{theorem}\label{cat-equiv}
  There is a homology equivalence $\Z^2 \times BO_{\infty,\infty}(\Z)
  \to \Omega BC$, and hence a homotopy equivalence of
  $\Omega^\infty$--spaces 
  \[
  \Z^2\times BO_{\infty,\infty}(\Z)^+ \simeq \Omega B\calC.
  \]
\end{theorem}
Here ``$+$'' denotes Quillen's plus construction with respect to the
commutator subgroup of $\pi_1 BO_{\infty,\infty}(\Z) \cong
O_{\infty,\infty}(\Z)$, which is perfect by \cite{Vaserstein} (or see
\cite[5.4.6, p.246]{Hahn-OMeara}).  The plus construction preserves
(generalized) (co)homology and kills a selected perfect subgroup of
$\pi_1$; in this case it kills the commutator subgroup and thus
abelianizes the fundamental group.  The space $\Z^2 \times
BO_{\infty,\infty}(\Z)^+$ is precisely the Hermitian algebraic
$K$--theory of the integers\footnote{Note that one must be careful
  whether one works with automorphisms of quadratic or symmetric
  bilinear forms; these groups are slightly different when 2 is not
  invertible in the ring.  However, it can be shown that the
  associated $K$--theories are rationally the same.}; it is known to
be an infinite loop space (see for example \cite{Loday-K-theory}).
Closely related to this, the Hermitian algebraic $K$--theory of
$\Z[1/2]$ has recently been studied by Berrick and Karoubi in
\cite{Berrick-Karoubi}; they compute the rational and 2--primary
homotopy groups.  The space $BO_{\infty,\infty}(\Z)^+$ is rationally
equivalent to $BO$, and it is equivalent to
$BO_{\infty,\infty}(\Z[1/2])^+$ away from the prime 2.

For the second variant of Morava's category, (slightly closer to the
original), let $\calC_{spin} \subset \calC$ denote the sub-2--category
with only \emph{spin} 4--manifolds (since they are simply connected,
this is equivalent to using only 4--manifolds with even intersection
form).  This spin sub-2--category is described in more detail in
section \ref{spin-section}.
\begin{theorem}\label{cat-equiv-spin}
  There is a map $\Z^2 \times B\Aut(\infty H \oplus \infty (-E_8)) \to
  \Omega B \calC_{spin}$ which is a homology equivalence away from the
  prime 2, and hence there is an $\Omega^\infty$--map
  \[
  \Z^2 \times B \Aut(\infty H \oplus \infty (-E_8))^+ \to \Omega B
  \calC_{spin}
  \]
  which is a homotopy equivalence away from the prime 2.
\end{theorem}
Here $E_8$ is the rank 8 irreducible form and $H$ is the standard rank 2
hyperbolic form 
\[
H = \left( \begin{array}{cc}
0 & 1 \\
1 & 0
\end{array} \right)
\]
Note that the commutator subgroup of $\Aut(\infty H
\oplus \infty (-E_8))$ is perfect---this follows from the argument in
\cite{Hahn-OMeara} given for $O_{\infty,\infty}(\Z)$.

The proofs of Theorems \ref{cat-equiv} and \ref{cat-equiv-spin} are
based on Tillmann's generalized group completion theorem
\cite{Tillmann-mcg-homotopy}, together with the isomorphisms of
Theorem \ref{stable-isomorphism} and Corollary \ref{twists-are-trivial2}.

From the perspective of physics most of the interesting mathematics is
related to the representation theory of the identity component of the
diffeomorphism group rather than the group of components.
Unfortunately our isotopy variants of the Morava categories lose sight
of this aspect entirely.  However, there is also an interest in
\emph{homotopy QFTs}, or \emph{flat QFTs}, which are essentially
representations of isotopy 2--categories such as ours.  See for
example \cite{Turaev}, \cite{Turner}.

The are operads closely related to the cobordism 2--categories we
define.  As an application of our analysis of mapping class groups in
dimension 4 we observe that Tillmann's higher genus surface operad
\cite{Tillmann-operad} has an a 4--manifold analogue.  Tillmann's
argument applies to this operad as well, so that the 4--manifold
mapping class group operad detects infinite loop spaces.

\subsection*{Organization of the paper}
The remainder of the paper is organized as follows.  In section
\ref{twist-section} we establish some properties of Dehn twists from
which we deduce Theorem \ref{twists-are-trivial1}.  In section
\ref{stable-groups} we review the results of Kreck, Wall and Quinn and
then combine them to give the proof of the the stable isomorphism
theorem, Theorem \ref{stable-isomorphism}.  The remainder of the paper
is concerned with the cobordism 2--categories of 4--manifolds.  In
section \ref{cats} we construct of the category $\calC$ and the map
into $K$--theory in detail, and in section \ref{cats-proof} study this
map and prove Theorem \ref{cat-equiv}.  Discussion of the spin case
and the proof of Theorem\ref{cat-equiv-spin} is contained in section
\ref{spin-section}.  In section \ref{operad-section} we discuss the
4--manifold mapping class group operad.

\section{Dehn twists}\label{twist-section}

The purpose of this section is to prove Theorem
\ref{twists-are-trivial1}, which gives a sufficient condition for when
the Dehn twist around a 3--sphere is actually isotopic (fixing the
boundary) to the identity.

By a \emph{Dehn twist} we shall mean the element of the mapping class
group constructed as follows.  The data required to construct a Dehn
twist is an embedding $S^3 \hookrightarrow M$ with trivial normal and
a loop $\alpha \in \Omega \Diff(S^3)$.  We think of the loop as
parametrized by the interval $(-1,1)$.  Let $V$ be a tubular
neighborhood of the embedded sphere.  One constructs a diffeomorphism
$\phi: M \to M$ by defining it to be the identity outside of $V$, and
on $V\cong S^3\times (-1,1)$ one sets $\phi(z,t) = \alpha_t(z)$.  One
easily sees that the isotopy class of $\phi$ depends only on the
homotopy class of $\alpha$ and the isotopy class of the embedding.

Let $M$ be a simply connected closed oriented 4--manifold , and let
$M'$ be the result of cutting out $n$ disjoint discs in $M$, so $M'$
is bounded by $n$ 3--spheres.  Our goal is to study the map $\Gamma(M')
\to \Gamma(M)$.  

\begin{proposition}\label{kernel-gen-by-dehn}
  The homomorphism $\Gamma(M') \to \Gamma(M)$ is surjective with
  kernel a quotient of $(\Z/2)^n$ generated by Dehn twists around the
  boundary components.
\end{proposition}

The proof will follow from a bit of elementary differential topology.
There is a fibration
\begin{equation}\label{eval-fib}
\Diff(M' \fix \partial M') \longrightarrow \Diff(M)
\stackrel{ev}{\longrightarrow} \Emb\left(\coprod^n D^4, M\right)
\end{equation}
where $\Emb$ is the space of embeddings (which extend to
diffeomorphisms of $M$).  For a single disc, linearization at the
center of the disc yields a homotopy equivalence between $\Emb(D^4,
M)$ and the frame bundle of $M$; see for instance \cite[Theorem
2.6.C]{Ivanov-survey}.  Similarly, when there is more than one disc
then there is a homotopy equivalence
\[ 
\Emb \left(\coprod^n D^4, M \right) \simeq FC_n(M)
\]
where $FC_n(M)$ is the \emph{framed configuration space}, consisting
of configurations of $n$ distinct ordered points in $M$ equipped with
(oriented) framings.  When $M$ is connected $FC_n(M)$ is also
connected.  Forgetting the framings gives a map to the usual
configuration space $C_n(M)$ of $n$ ordered points which fits into a
fibration
\begin{equation}\label{framed-unframed-fib}
SO(4)^n \to FC_n(M) \to C_n(M).
\end{equation}

\begin{lemma}
  If $M$ is closed and simply connected then $\pi_1 FC_n(M)$ is a
  quotient of $(\Z/2)^n$, with the generators corresponding to
  rotations of each of the framings.
\end{lemma}
\begin{proof}
  Since $M$ is 4--dimensional, removal of a finite number of points in
  $M$ preserves simple-connectedness and connectedness.  Induction on
  $k$ with the Faddell--Neuwirth fibrations $M - (\mbox{$k$ points})
  \to C_{k+1}(M) \to C_{k}(M)$ shows that the unframed configuration
  spaces are all simply connected.  The result then follows from the
  homotopy exact sequence of the fibration (\ref{framed-unframed-fib}).
\end{proof}

From the fibration (\ref{eval-fib}) there is an exact sequence of
homotopy groups,
\begin{equation}
  \pi_1 \Diff(M) \stackrel{\alpha}{\to} \pi_1 FC_n(M) \stackrel{\delta}{\to} 
  \Gamma(M') \to \Gamma(M) \to 0.
\end{equation}
Note that a rotation of a framing in $\pi_1 FC_n(M)$ is sent by
$\delta$ to the Dehn twist around the corresponding boundary sphere of
$M'$ in $\Gamma(M')$, and Proposition \ref{kernel-gen-by-dehn} follows.
From this we also see that,
\begin{lemma}
  The map $\Gamma(M') \to \Gamma(M)$ is an isomorphism if and only if
  $\alpha$ is surjective.
\end{lemma}

Fixing distinct points $p_1, \ldots p_n \in M$, the diagram of
fibrations
\begin{equation*}
\begin{diagram}
\node{\Diff\left(M \fix \bigcup p_i \right)} \arrow{e} \arrow{s} 
\node{\Diff(M)} \arrow{e} \arrow{s} \node{C_n(M)} \arrow{s} \\
\node{SO(4)^n } \arrow{e} \node{FC_n(M)} \arrow{e} \node{C_n(M)}
\end{diagram}
\end{equation*}
induces a homomorphism of long exact sequences, 
\begin{equation*}
\begin{diagram}
\node{\pi_2 C_n(M)} \arrow{s,r}{id} \arrow{e}
\node{\pi_1 \Diff\left(M \fix \bigcup p_i\right)} \arrow{e} \arrow{s,r}{\beta}
\node{\pi_1 \Diff(M)} \arrow{s,r}{\alpha} \arrow{e} \node{0}\\
\node{\pi_2 C_n(M)} \arrow{e} \node{\pi_1 SO(4)^n}
\arrow{e} \node{\pi_1 FC_n(M)} \arrow{e} \node{0} 
\end{diagram}
\end{equation*}
from which we see that,
\begin{lemma}\label{surjectivity}
$\alpha$ is surjective if and only if $\beta$ is.
\end{lemma}
\begin{proof}
One direction is immediate and the other follows from the Five Lemma.
\end{proof}

\begin{lemma}\label{CPlemma}
  $\Gamma(\C P^2 - \{\mbox{2 discs}\}) \cong \Gamma(\C P^2)$.  
\end{lemma}
\begin{proof}
  By Lemma \ref{surjectivity} is suffices to show that $\beta$ is
  surjective.  We do this by constructing $S^1$--actions on $\C P^2$
  which hit each of the generators of $\Z/2\times \Z/2
  \twoheadrightarrow \pi_1 FC_2(\C P^2)$.  Let $p_1 = [0, 0, 1] \in \C
  P^2$ and $p_2 = [1, 0, 0] \in \C P^2$, and consider the $S^1$--action
  defined by $\lambda \cdot [x, y, z] = [x, y, \lambda z]$.  This
  action fixes both $p_1$ and $p_2$, and hence there is a
  representation of $S^1$ on the tangent space at each of these two
  points.  The complex 1--dimensional representations of $S^1$ are
  labelled by the integers and one can easily identify the
  representations on the tangent spaces of the fixed points by
  choosing local coordinates: the representation on $T_{p_1}\C P^2$ is
  $({\bf -1})\oplus ({\bf -1})$, and at $p_2$ the representation is
  $({\bf 0}) \oplus ({\bf 1})$.  These correspond to compositions of
  group homomorphisms
  \[
    S^1 \stackrel{\rho_i}{\longrightarrow} U(1) \times U(1)
    \stackrel{\iota}{\hookrightarrow} U(2) \hookrightarrow SO(4).
  \]
  The induced maps on fundamental groups are
  \[
  \Z \stackrel{\rho_{i*}}{\longrightarrow} \Z\times \Z 
  \stackrel{\iota_*}{\longrightarrow} \Z \twoheadrightarrow \Z/2
  \]
  with $\iota_*$ being addition and $\rho_{i*}$
  determined by the representation of $S^1$.  That is, 
  \begin{eqnarray*}
    \rho_{1*} : n \mapsto (-n, -n),\\
    \rho_{2*} : n \mapsto (0,n).
  \end{eqnarray*}
  Thus $(0,1)$ is in the image of
  \[
  \beta: \pi_1 \Diff\left(\C P^2 \fix \bigcup p_i\right)
  \longrightarrow \Z/2\times\Z/2 \twoheadrightarrow \pi_1 FC_2(\C
  P^2).
  \]
  By letting $S^1$ act instead on the first coordinate of $\C P^2$,
  one sees that $(1,0)$ is also in the image of $\beta$.
\end{proof}

As a consequence we have that a Dehn twist around either of the
boundary components in $\C P^2 - \{\mbox{2 discs}\}$ is isotopic
(keeping the boundary fixed) to the identity.

The spin case in Theorem \ref{twists-are-trivial2} is handled by the
following lemma.  

\begin{figure}\caption{\label{1-skel} The relative $1$-skeleton of $(M',\partial M')$.}
\begin{center}
\vspace{5mm}
\input{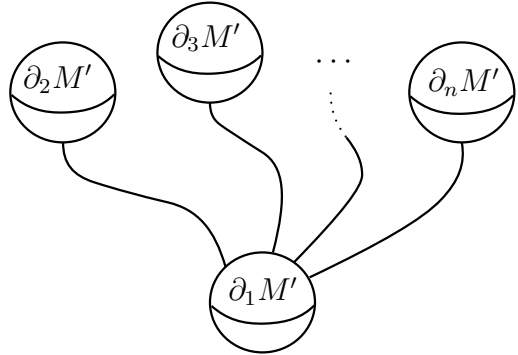}
\end{center}
\end{figure}

\begin{lemma}\label{spin-dehn2}
Suppose $M$ is a smooth closed oriented simply connected 4--manifold
which is spin, and $M'$ is the result of removing the interiors of $n$
disjoint discs.  There is a homomorphism $\Gamma(M') \to
H^1(M',\partial M'; \Z/2) \cong (\Z/2)^{n-1}$ whose restriction to the
kernel $K$ of $\Gamma(M') \to \Gamma(M)$ is surjective.
\end{lemma}
\begin{proof}
Let $s$ be the unique spin structure on $M$.  Spin structures on $M'$
(relative to the boundary) are an affine space over $H^1(M',\partial
M';\Z/2)$, and there is a natural choice of basepoint given by the
restriction of $s$ to $M'$.  Now, define a homomorphism $\Gamma(M')
\to H^1(M',\partial M';\Z/2)$ by
\[
\varphi \mapsto \varphi^*(s) - s \in H^1(M',\partial M';\Z/2);
\]
this measures the extent to which $\varphi$ fails to preserve the spin
structure $s$.  by Proposition \ref{kernel-gen-by-dehn} $K$ is
generated by Dehn twists around the boundary spheres, and these twists
act transitively on the set of spin structures on $(M',\partial M')$.
To see this, think of the relative 1--skeleton of $M'$ as a union of
arcs joining one fixed boundary sphere with each of the other boundary
spheres as in Figure \ref{1-skel}. Spin structures on $(M',\partial
M')$ correspond to labelings of the arcs in the relative 1--skeleton
by elements of $\Z/2$.  If $\varphi$ is a
Dehn twist around a boundary sphere $S$ then $\varphi^*$ reverses the
label of each arc having an endpoint on $S$.  Hence the homomorphism
defined above is surjective.
\end{proof}

\begin{proof}[Proof of Theorem \ref{twists-are-trivial1}]
  Recall the setup: $M$ is an arbitrary 4-manifold containing $\C P^2$
  as a connected summand.  Let $S \hookrightarrow M$ be a 3--sphere
  embedded in $M$ with trivial normal bundle and let $\alpha$ denote
  the Dehn twist around $S$.  We cut $M$ along $S$, producing a
  manifold $M'$, and we now regard $\alpha$ as Dehn twist around one
  of the new boundary components.  Up to diffeomorphism, we may assume
  that the boundary component around which $\alpha$ twists lies on the
  $\C P^2$ summand, as in figure \ref{untwist-fig}, so $\alpha$ is in
  the image of the composition
  \[
  \Gamma(\C P^2 - \{\mbox{2 discs}\}) \to \Gamma(M') \to \Gamma(M).
  \]
  Lemma \ref{CPlemma} tells us that the element of $\Gamma(\C P^2 -
  \{\mbox{2 discs}\})$ which maps to $\alpha$ is the zero element, so
  $\alpha$ is zero in $\Gamma(M)$.  The same argument holds with $\C
  P^2$ replaced by $\overline{\C P^2}$ throughout.
\end{proof}

\begin{proof}[Proof of Corollary \ref{twists-are-trivial2}]
  We have $M'$ obtained from $M$ by removing a collection of discs.
  The map $\Gamma(M') \to \Gamma(M)$ given by gluing the discs back in
  is surjective with kernel $K$.  By Proposition
  \ref{kernel-gen-by-dehn} $K$ is generated entirely of Dehn twists
  around boundary 3--spheres.  If $M$ contains $\C P^2$ or
  $\overline{\C P^2}$ as a connected summand then all Dehn twists are
  isotopic to the identity by Theorem \ref{twists-are-trivial1}, so
  $K=0$.  If $M$ is spin then by Lemma \ref{spin-dehn2} $K$ is
  isomorphic to either $(Z/2)^n$ or $(\Z/2)^{n-1}$ since it is a
  quotient of $(\Z/2)^n$.
\end{proof}

\begin{figure}\caption{\label{untwist-fig}}
\begin{center}
\vspace{5mm}
\input{untwist-new.pst}
\end{center}
\end{figure}

In the spin case of Corollary \ref{twists-are-trivial2} it appears to
be a difficult problem to decide in general precisely how many powers
factors of $\Z/2$ there are in $K$.  As easy case is conencted sums of
$S^2\times S^2$.

\begin{proposition}\label{spin-dehn1}
If $M$ is a connected sum of copies of $S^2\times S^2$ and $M'$ is
obtained from $M$ by deleting the interiors of $n$ discs then
$\mathrm{ker}(\Gamma(M') \to \Gamma(M)) = K \cong \Z/2^{n-1}$.
\end{proposition}
\begin{proof}
First consider the homomorphism
\begin{equation}\label{two-spheres-one-disc}
\Gamma(S^2 \times S^2 - \{\mbox{a disc}\}) \to \Gamma(S^2 \times S^2)
\end{equation}
This is actually an isomorphism.  The idea of the proof is the same as
for Lemma \ref{CPlemma}; we look for a circle action on $S^2\times
S^2$ which fixes a point $p$ and rotates the tangent space $T_p$
through the nontrivial element of $\pi_1 SO(4)$.  Rotation about the
polar axis on the first sphere, with $p=(\mbox{north pole},\mbox{north
pole})$ does the job.

From the isomorphism (\ref{two-spheres-one-disc}) it follows that if
$M$ is a connected sum of copies of $S^2 \times S^2$ then gluing in a
disc on $M'$ kills one factor of $\Z/2$ in $K$ as long as there is at
least one remaining boundary component because the Dehn twist around
any boundary sphere is isotopic to the product of the twists around
each of the other boundary spheres.  The conclusion now follows by
induction on $n$.
\end{proof}

A connected sum of copies of $S^2\times S^2$ is essentially the only
case which can be handled directly by circle actions in light of the
classification of locally smooth circle actions on closed simply
connected 4--manifolds given in \cite{Fintushel}. If such an action
exists then the manifold must be a connected sum of copies of $S^2
\times S^2$, $\C P^2$, $\overline{\C P^2}$, or a homotopy 4--sphere.

\section{The stable groups  $\Aut_\infty$ and $\Gamma_\infty$}
\label{stable-groups}
In this section we recall the theorems of Wall, Kreck, and Quinn and
show how they combine to give a lifting of $\Aut(Q_M)$ into
$\Gamma(MX^k)$ which stabilizes to an inverse of $\Gamma_\infty(M) \to
\Aut_\infty(Q_M)$.

In the statements of the following three theorems $M$ shall be a
simply connected compact oriented smooth 4--manifold, possibly bounding
some number of homology 3--spheres.

\begin{theorem}{\cite{Wall2}}\label{wallthm}
  If $M$ is of the form $N\# (S^2 \times S^2)$ with $Q_N$ either
  indefinite or of rank $\leq 8$ then $\Gamma(M) \to \Aut(Q_M)$ is
  surjective.
\end{theorem}

Two diffeomorphisms $\varphi_0,\varphi_1$ of $M$ are said to be
\emph{pseudo-isotopic} if they are the restrictions to $0\times M$ and
$1\times M$ respectively of a diffeomorphism $\Phi$ of $I\times M$.
Pseudo-isotopy is in general a coarser equivalence relation than
isotopy, but is finer than homotopy.  Let $P(M)$ denote the group of
pseudo-isotopy classes of diffeomorphisms of $M$; this group is a
quotient of $\Gamma(M)$, and the morphism from $\Gamma(M) $to
$\Aut(M)$ descends to $P(M)$.

\begin{theorem}{\cite{Kreck}}\label{kreckthm}
  $P(M) \to \Aut(Q_M)$ is always injective, provided that the
  pseudo-isotopies are \emph{not} required to fix the boundary
  point-wise.
\end{theorem}

Thus $P(M) \cong \Aut(Q_M)$ whenever $M$ satisfies the conditions of
Wall's theorem and one is lead to ask what the relationship between
$\Gamma(M)$ and $P(M)$ is.  Quinn has given a good answer to this
question; he proves that in dimension 4 pseudo-isotopy implies
isotopy-after-stabilization.

\begin{theorem}{\cite{Quinn}}\label{quinnthm}
  If $\varphi \in \Diff(M)$ is pseudo-isotopic to the identity then
  for $k$ large enough its extension (by identity on $S^2\times S^2$)
  to an element of $\Diff(M\#k(S^2\times S^2))$ is isotopic to the
  identity.  (All isotopies and pseudo-isotopies fix the boundary of
  $M$ pointwise.)
\end{theorem}
\begin{remark}
  For our purposes we require only this weak version of Quinn's
  theorem---the stronger form actually gives existence of an isotopy
  not merely between the $t=1$ end of the pseudo-isotopy and the
  identity diffeomorphism on $M$, but between the pseudo-isotopy
  (regarded as an element of $\Diff(M\#k(S^2 \times S^2)\times I,
  \partial M\ \times I)$) and the identity on $M\#k(S^2\times S^2)
  \times I$.
\end{remark}

Note that the value of $k$ may depend on $\varphi$; it could
potentially be unbounded as $\varphi$ ranges over all pseudo-isotopy
classes.  However, in all known examples $k=1$ suffices.

Quinn's theorem applies even when $M$ has boundary and all
diffeomorphisms and (pseudo)-isotopies are taken to fix the boundary
point-wise, as does Wall's surjectivity result; however, the
injectivity result of Kreck is no longer valid if one requires the
boundary to be fixed because there can be a kernel consisting of
twists around boundary components which induce isomorphism on
homology. This limitation is precisely what necessitates our analysis
of Dehn twists in the previous section since our stabilization process
requires that diffeomorphisms fix the boundary point-wise in order to
have a well-defined extension.

We now consider stabilization of $\Gamma_M$ and $\Aut(Q_M)$ by gluing
on copies of $X=\C P^2 \# \overline{\C P^2} - \{\mbox{2 discs}\}$.

\begin{lemma}\label{maintool}
  Suppose $M$ is a simply connected oriented smooth compact
  4--manifold bounded by a collection of 3--spheres.  Then each
  automorphism of $Q_M$ is induced (as an element of $\Aut_\infty(Q_M)$)
  by a unique element of $\Gamma_\infty(M)$.
\end{lemma}
\begin{proof}
  Let $\alpha \in \Aut(Q_M)$.  Even if $\alpha$ is not induced by a
  diffeomorphism, Theorem \ref{wallthm} implies that for $\ell \geq 2$
  its image in $\Aut(Q_{MX^\ell})$ is induced by a diffeomorphism.
  Any two diffeomorphisms representing $\alpha$ are pseudo-isotopic
  after closing off the boundary with discs by Theorem \ref{kreckthm}.
  Furthermore, by Theorem \ref{quinnthm} these two representatives are
  actually isotopic after extending to $MX^k$, for some $k\geq \ell$
  large enough, and then closing off the boundary with discs.  But the
  extensions are isotopic even before closing off the boundary by
  Corollary \ref{twists-are-trivial2}.
\end{proof}
This next lemma follows immediately from the previous.
\begin{lemma}\label{maintool2}
  There is a unique inclusion $\Aut(Q_M)
  \hookrightarrow \Gamma_\infty(M)$ such that the following
  diagram commutes,
\begin{equation*}
\begin{diagram}
\node[2]{\Gamma_\infty(M)} \arrow{s} \\
\node{\Aut\left(Q_M\right)} \arrow{ne,J} \arrow{e,J}
\node{\Aut_\infty(Q_M).}
\end{diagram}
\end{equation*}
\end{lemma}

\begin{proof}[Proof of Theorem \ref{stable-isomorphism}]
  We produce an inverse to the map $\Gamma_\infty(M)\to
  \Aut_\infty(Q_M)$ by exhibiting compatible maps on each of the terms
  in the directed system which defines $\Aut_\infty(Q_M)$.  The liftings
  of Lemma \ref{maintool2} serve this purpose.  We need only check that
  the diagrams
\begin{equation*}
\begin{diagram}
\node{\Aut(Q_{MX^k})} \arrow{s,J} \arrow{se,J} \\
\node{\Aut(Q_{MX^{k+1}})} \arrow{e,J} \node{\Gamma_\infty(M)}
\end{diagram}
\end{equation*}
all commute, but this follows immediately from Lemmas \ref{maintool}
and \ref{maintool2}.  Hence we obtain a homomorphism
$\Aut_\infty(Q_M)\to \Gamma_\infty(M)$ which is the desired inverse by
construction. 

We now claim that the stabilized automorphism group $\Aut_\infty(Q_M)$
is isomorphic (non-canonically) to the split-signature orthogonal
group $O_{\infty,\infty}(\Z) = \colim_k O_{k,k}(\Z)$, which is clearly
independent of the initial manifold $M$.  By the Hasse-Minkowski
classification of odd indefinite unimodular forms, $Q_M \oplus (1)
\oplus (-1) \cong a(1) \oplus b(-1)$ for some natural numbers $a,b$.
Thus
\begin{align*} 
\Aut_\infty(Q_M) & = \Aut_\infty(Q_M\oplus (1) \oplus (-1)) \\
                 & \cong \colim_{k\to \infty} O_{a+k,b+k}(\Z) \\
                 & \cong O_{\infty,\infty}(\Z).
\end{align*}
Note that passing from the first to the second line requires choosing
a basis of vectors of length 1, so the resulting identification is
non-canonical.  In passing to the third line we have used the fact
that $\colim_k O_{k,k}(\Z) \cong \colim_k O_{k,k+n}(\Z)$; this can be
proved by comparing both sides to $\colim_{j,k} O_{j,k}(\Z)$.  The two
embeddings $\N \hookrightarrow \N\times \N$ given by $k\mapsto (k,k)$
and $k \mapsto (k,k+n)$ are both co-final subsystems.
\end{proof}
 
In applying Tillmann's group completion theorem to identify the
homotopy type of the cobordism categories described in the next
section we will need an additional result about this group which we
record now.  Let $L$ be a simply connected smooth compact oriented
4--manifold with outgoing boundary (a union of 3--spheres) compatible
with the incoming boundary of $M$ so that $LM$ is a well defined
composition.
\begin{lemma}\label{translation-iso}
The inclusion $i: \Aut_\infty(Q_M) \hookrightarrow \Aut_\infty(Q_{LM})
= \Aut_\infty(Q_L\oplus Q_M)$ is an integral homology equivalence.
\end{lemma}
\begin{proof}
The Hasse-Minkowski classification implies that 
\begin{align*}
Q_{LMX^2} & \cong (Q_L \oplus Q_X) \oplus (Q_M \oplus Q_X) \\
          & \cong (\ell_1 (1) \oplus \ell_2 (-1)) \oplus (m_1 (1) \oplus m_2(-1)),
\end{align*}
so it suffices to show that \[ \colim_{k\to \infty} O_{k,k+n}(\Z) \hookrightarrow
\colim_{k \to \infty} O_{k+\ell_1, k+n+\ell_2}(\Z)\] is a homology isomorphism.

The group $O_{k+\ell_1,k+n+\ell_2}$ contains $\Sigma_{k+\ell_1} \times
\Sigma_{k+n+\ell_2}$ as permutations of basis vectors.  Thus for some
$k'$ large enough there exists a basis permutation which conjugates
the image of the stabilization embedding $O_{k+\ell_1,k+n+\ell_2}(\Z)
\hookrightarrow O_{k' +\ell_1, k'+n + \ell_2}(\Z)$ into the image of
the embedding $O_{k',k'+n}(\Z) \hookrightarrow O_{k'+\ell_1,
k'+n+\ell_2}(\Z)$.  Since conjugation induces an isomorphism on group
homology, and each class in $H_*O_{\infty,\infty}(\Z)$ comes from
$H_*O_{k,k+n}(\Z)$ for some $k$ large enough, this proves
surjectivity.  Injectivity follows from a similar argument.  Note that
this argument can be easily modified to work also in the case of even
indefinite forms.
\end{proof}

Since $\Aut_\infty(Q_M) \cong \Gamma_\infty(M)$, we also have:
\begin{corollary}\label{translation-iso2}
  $\Gamma_\infty(M) \to \Gamma_\infty(LM)$ is also an integral
  homology equivalence.
\end{corollary}

Note that the analogue of Corollary \ref{translation-iso2} for
surfaces holds by virtue of Harer stability.  In dimension 4 we are
able to replace the need for homological stability with general
properties of the homology of linear groups.

\section{Definitions of the categories of 4--manifolds}\label{cats}

In this section we construct the 2--category $\calC$ and the map into
$K$-theory which comes from a 2--functor into a closely related
category $\calK$.  The definitions we use are natural extensions of
the definitions found in \cite{Tillmann-mcg-homotopy} and
\cite{Tillmann-splitting}.  In particular, our 2--category $\calC$ is
constructed precisely along the lines of the surface cobordism
2--category in the second reference above.

Both $\calC$ and $\calK$ have the same underlying ordinary category
(i.e. the same objects and morphisms) which we denote by $\calC_0$;
conceptually, this category should be thought of as the cobordism
category of (unions of) 3--spheres and simply connected oriented
4--manifolds.  However, one must be careful in defining the morphisms
so that composition is well-defined and the result is a small
category.  The 2--morphisms of $\calC$ and $\calK$ will be isotopy
classes of diffeomorphisms and isomorphisms of the intersection form,
respectively.  Let us proceed in detail.

\vspace{0.5cm}
\emph{Objects of $\calC_0$:} The objects are the non-negative integers,
with $n$ thought of as representing a disjoint union of $n$ copies of
$S^3$.

\vspace{0.5cm} \emph{Morphisms of $\calC_0$:} Let $\mathscr{A}(n)$
denote a set of manifolds containing precisely one representative from
each diffeomorphism class of compact oriented connected
simply connected 4--manifolds bounded by $n+1$ ordered 3--spheres.  We
equip each boundary sphere with a collar and we consider the first $n$
boundary components as \emph{in-going} and the final boundary
component as \emph{out-going}.

We now allow these \emph{atomic} manifolds to freely generate the
morphism sets via finite sequences of the three operations of:
\begin{itemize} 
\item[\bf 1] Gluing the out-going boundary of one morphism to the
  in-going boundary of another using the ordering of the boundary
  components and the collars. 
\item[\bf 2] Taking disjoint unions.
\item[\bf 3] Renumbering the in-going and out-going boundary
  components.
\end{itemize}
The morphism set $\calC_0(m,n)$ consists of all such composites with
$m$ incoming and $n$ outgoing boundary components respectively,
together with an identity morphism when $m=n$.  Composition of
morphisms is given by gluing the out-going end of one to the in-going
end of the other.  Given $M\in \calC_0(a,b)$ and $N\in \calC_0(c,d)$, the
ordering of the in-going boundary components of a disjoint union $M
\sqcup N \in \calC_0(a+c, b+d)$ is determined by taking first the
in-going boundary components of $M$ followed by those of $N$, and
likewise for the out-going boundary ordering. 

\begin{remark}
  We have imposed the requirement that each component has precisely
  one outgoing boundary sphere to ensure that all compositions remain
  simply connected.
\end{remark}

It may at first seem strange that the atomic manifolds are allowed to
\emph{freely} generate the morphisms, since clearly a given cobordism
$M^4$ can be written as a composition of smaller pieces in many
different ways and one would want the different decompositions of $M$
to all represent the same cobordism.  However, we follow an aspect of
the philosophy of 2-categories; rather than trying to \emph{equate}
all of the different decompositions of $M$ into atomic manifolds, we
use the 2-morphisms to encode the property that the different
decompositions are all isomorphic.

\begin{definition}
\label{CDef}
Let $\calC$ be the (strict) 2--category with underlying category
$\calC_0$, and with 2--morphisms given by isotopy classes (using
isotopies constant on the boundary) of diffeomorphisms which respect
the parametrizations and ordering of the boundary components.  Thus if
$M$ and $N$ are morphisms that are diffeomorphic (respecting the
boundary data), then the 2--morphisms from $M$ to $N$ are the isotopy
classes of diffeomorphisms from $M$ to $N$ which respect the boundary
data; if $M$ and $N$ are not diffeomorphic then the set of
2--morphisms between them is empty.  In the case where $M$ is a
morphism that is obtained from an identity morphism by renumbering the
boundary we take the 2--morphisms from $M$ to $N$ to be empty unless
$N=M$, in which case there is only the identity 2--morphism.
Horizontal composition of 2-morphisms is induced by gluing; vertical
composition of 2--morphisms is induced from composition of
diffeomorphisms.
\end{definition}

\begin{definition}
\label{Kdef}
The (strict) 2--category $\calK$ has underlying category $\calC_0$ and
the 2--morphisms are now the isomorphisms of the intersection forms: a
2--morphism from $M$ to $N$ is an isomorphism $H_2 M
\stackrel{\cong}{\longrightarrow} H_2 N$ which preserves the
intersection form.
\end{definition}

One may form simplicial categories $\mathscr{B}\calC$ and
$\mathscr{B}\calK$ by replacing the morphism categories in $\calC$ and
$\calK$ with their nerves.  The nerve of a simplicial category is a
bisimplicial set; the geometric realization of a simplicial category
is defined to be the geometric realization of this bisimplicial nerve.  
We thus define the geometric realization $B\calC$
($B\calK$) of the 2--category $\calC$ ($\calK$ resp.) as the geometric
realization of the associated simplicial category;
\[
B\calC := B(\mathscr{B}\calC), \mbox{         }
B\calK := B(\mathscr{B}\calK).
\]

Disjoint union provides a strict symmetric monoidal product on each of
$\calC$ and $\calK$ and hence an infinite loop structure on their
geometric realizations (since both spaces are connected).  There is an
obvious natural 2--functor $F: \calC \to \calK$; it is identity on
objects and morphisms, and it sends the isotopy class of a
diffeomorphism $\phi: M {\to} N$ to the induced isomorphism of
intersection forms $\phi_*: (H_2 M, Q_M) {\to} (H_2 N, Q_N)$.

\section{The proof of Theorem \ref{cat-equiv}}\label{cats-proof}

In this section we prove Theorem \ref{cat-equiv} by identifying the
homotopy type of $\Omega B\calK$ and showing that the 2--functor $\calC
\to \calK$ is a homotopy equivalence after group completion.  The
proof is based on a group-completion argument, closely following
Tillmann's proof \cite{Tillmann-mcg-homotopy} that $\Z \times
B\Gamma_\infty^+$ is an infinite loop space.  At the conceptual level,
where Tillmann uses Harer--Ivanov stability of mapping class groups to
obtain a homology fibration we substitute Lemma \ref{translation-iso}
together with Theorem \ref{stable-isomorphism}.

Fix a morphism $X \in \calC_0(1,1)$ diffeomorphic to $\C P^2 \#
\overline{\C P^2}$ with two discs removed, and consider the
contravariant functor $\mathcal{X}_\infty: \calC_0 \to
(\mbox{Simplicial Sets})$ defined by the telescope construction:
\[
  \mathcal{X}_\infty(n) := \mathrm{hoColim} \{ B\calC(n,1)
  \stackrel{X}{\longrightarrow} B\calC(n,1)
  \stackrel{X}{\longrightarrow} \cdots \},
\]
where $B\calC(n,1)$ is the simplicial nerve of the morphism category
$\calC(n,1)$.
A morphism $L:m \to n$ induces a map
$\mathcal{X}_\infty(n) \to \mathcal{X}_\infty(m)$ by gluing on the
left.  Similarly we set
\[
  \mathcal{Y}_\infty(n) := \mathrm{hoColim} \{ B\calK(n,1)
  \stackrel{X}{\longrightarrow} B\calK(n,1)
  \stackrel{X}{\longrightarrow} \cdots \}.
\]

\begin{lemma}\label{telescope}
$\mathcal{Y}_\infty(n) \simeq \Z^2 \times BO_{\infty,\infty}(\Z)$
\end{lemma}
\begin{proof}
  Swapping the order of the classifying space functor and the homotopy
  colimit functor expresses $\mathcal{Y}_\infty(n)$ as the classifying
  space of the colimit of the hom-set groupoids $\calK(n,1)$.  
  The connected components of the groupoid $\calK(n,1)$ correspond the
  set of isomorphism classes of intersection forms (of any rank and
  signature).  Stabilization by gluing on copies of $X$ corresponds to
  block addition with the intersection form $(1)\oplus(-1)$, and under
  this stabilization two intersection forms eventually become
  isomorphic if and only if they have the same rank and the same signature.
  Thus the connected components of the colimit groupoid 
  \[
  \colim \{\calK(n,1) \to \calK(n,1) \to \cdots \}
  \]
  are $\Z^2$, the group completion of the additive monoid formed by
  the pairs $\{(\mbox{rank}, \mbox{signature})\} \subset \N \times \Z$.
  Hence the components of $\mathcal{Y}_\infty(n)$ are in bijection
  with $\Z^2$ and one sees that each component of
  $\mathcal{Y}_\infty(n)$ is the classifying space of a groupoid with
  underlying group $O_{\infty,\infty}(\Z)$.
\end{proof}

\begin{lemma}\label{telescope-equiv}
  $\mathcal{X}_\infty(n) \simeq \Z^2 \times BO_{\infty,\infty}(\Z)$
  and the 2--functor $\calC \to \calK$ induces a homotopy equivalence
  $\mathcal{X}_\infty(n) \simeq \mathcal{Y}_\infty(n)$.
\end{lemma}
\begin{proof}
  In $\mathcal{X}_\infty$ swap the classifying space functor with the
  homotopy colimit functor, so 
  \[
  \mathcal{X}_\infty(n) = B \colim \{ \calC(n,1) \to \calC(n,1) \to
  \cdots\}.
  \]
  The connected components of the groupoid $\calC(n,1)$ correspond to
  the diffeomorphism classes of cobordisms $n \to 1$ and
  are thus indexed by the atomic manifolds $\mathscr{A}(n,1)$.  As we
  stabilize by gluing on copies of $X$, two objects in the groupoid
  $\calC(n,1)$ eventually become isomorphic if and only if their
  underlying cobordisms eventually become diffeomorphic.  By Wall's
  stable diffeomorphism classification \cite{Wall2}, this happens if
  and only if the two objects have intersection forms which are stably
  isomorphic.  Hence, as in Lemma \ref{telescope}, the set of
  components of $\mathcal{X}_\infty(n)$ is $\Z^2$, corresponding to
  the rank and signature.  Each component of $\mathcal{X}_\infty(n)$
  is easily seen to be the classifying space of the stable 4--manifolds
  mapping class group, $B\Gamma_\infty$---here we use
  \ref{stable-isomorphism} to know that this group is independent of
  the initial manifold and hence all components are homotopy
  equivalent.

  Now the 2--functor $\calC \to \calK$ clearly induces a bijection on
  components, and on each component it induces the natural map
  $B\Gamma_\infty \to B\Aut_\infty$ which is a homotopy equivalence by
  Theorem \ref{stable-isomorphism}.
\end{proof}

The functor $\mathcal{X}_\infty$ ($\mathcal{Y}_\infty(n)$) is a
$B\calC$--diagram ($B\calK$--diagram, respectively) in the language of
\cite{Tillmann-mcg-homotopy}.  That is to say, the simplicial set
$$\coprod_n \mathcal{X}_\infty(n)$$ is equipped with a unital and associative
simplicial action of $B\calC$:
$$B\calC(n,m)\times \mathcal{X}_\infty(m) \to \mathcal{X}_\infty(n)$$
defined by composition on the left.  One may thus form the simplicial
Borel construction (a.k.a the bar construction)
\[
  (E_{B\calC}\mathcal{X}_\infty)_k =
  \underbrace{B\calC(-,-)\times_{\N}\cdots \times_{\N} B\calC(-,-)}_{k
  \mbox{ factors}} \times_{\N}\mathcal{X}_\infty(-).
\] 
The Borel construction commutes with the telescope, so
\[
  E_{B\calC} \mathcal{X}_{\infty} = \operatorname{hoColim}_X \{ E_{B\calC}
  B\calC(-,1)\}.
\]
As observed by Tillmann, $E_{B\calC} B\calC(-,1)$ is precisely
the nerve of the comma category $(B\calC \downarrow 1)$ of objects in
$B\calC$ over 1.  This latter category is contractible because it
contains the identity $1\to 1$ as a terminal object, and hence
$E_{B\C}\mathcal{X}_\infty$ is contractible as it is a homotopy
colimit (over a contractible category) of contractible spaces.

For each $n$ we have a pull-back diagram
\begin{equation*}
\begin{diagram}
  \node{\mathcal{X}_\infty(n)} \arrow{e} \arrow{s}
  \node{E_{B\calC}\mathcal{X}_\infty} \arrow{s} \\
  \node{n} \arrow{e,J} \node{B\calC}
\end{diagram}
\end{equation*}
and thus there is a map into the homotopy fibre: 
\begin{equation}\label{group-completion}
\mathcal{X}_\infty(n) \to \Omega B\calC.
\end{equation} 
The left translation maps $L\circ: \mathcal{X}_\infty(n) \to
\mathcal{X}_\infty(m)$ are all integral homology equivalences by
Theorem \ref{stable-isomorphism} and Lemma \ref{translation-iso}, so
Tillmann's generalized group completion theorem implies that the group
completion map (\ref{group-completion}) is an integral homology
equivalence.

Replacing $\mathcal{X}$ with $\mathcal{Y}$ and $\calC$ with $\calK$,
we obtain a homology equivalence
\[
\mathcal{Y}_\infty(n) \to \Omega B \calK.
\]
The homology equivalence of Theorem \ref{cat-equiv} now follows from
Lemma \ref{telescope-equiv}, and the homotopy equivalence then follows
from the properties of the plus construction together with the Whitehead
theorem for simple spaces.

We note here that the infinite loop space structure on $\Omega B
\calK$ coming from the symmetric monoidal product in $\calK$ coincides
with the usual infinite loop structure induced from direct sum.

\section{The spin case}\label{spin-section}

To obtain the proof of Theorem \ref{cat-equiv-spin} we simply restrict
everything in sight to even intersection forms and then check that the
proof of Theorem \ref{cat-equiv} goes through, at least away from the
prime 2.

Let us proceed in more detail. Thus $\calC_{spin}$ is the
sub-2-category of $\calC$ with the same objects and containing only
those morphisms which have even intersection forms; since all
4--manifolds here are simply connected this is exactly equivalent to
admitting a spin structure, and such a structure is unique whenever it
exists.  The 2--morphisms are all diffeomorphisms whose source and
target are even; since we are dealing with simply connected spin
manifolds, any diffeomorphism automatically respects the spin
structure.  Similarly, $\calK_{spin}$ is the restriction of $\calK$ to
even morphisms.

We would like to form spin analogues of the telescopes
$\mathcal{X}_\infty$ and $\mathcal{Y}_\infty$, but $\CPCPB$ is no
longer available in the spin setting since it has odd intersection
form, so instead we stabilize using $S^2\times S^2$.  For spin mapping
class groups the analogue of Corollary \ref{translation-iso2} is now
only a homology isomorphism with $\Z[1/2]$ coeffficients.  Corollary
\ref{twists-are-trivial2} no longer provides an isomorphism, so in the
proof of Lemma \ref{maintool} there is an indeterminacy when passing
from an element in the mapping class group of a closed manifold to an
element in the mapping class group of that manifold minus some discs.
However, by Lemma \ref{kernel-gen-by-dehn}, this indeterminacy is
purely 2--torsion.  Letting $\Gamma_\infty^{spin}(M)$ and
$\Aut_\infty^{spin}(Q_M)$ denote the mapping class group and
automorphism group stabilized with $S^2\times S^2$, we therefore have,
\begin{theorem}
The map $\Gamma_\infty^{spin}(M) \to \Aut_\infty^{spin}(Q_M)$ is
surjective with kernel of exponent 2.  More precisely, the kernel is
either $(\Z/2)^{n-1}$ or $(\Z/2)^{n}$, where $n$ is the number of
boundary spheres that $M$ has (one of these is used for
stabilization).
\end{theorem}

Note that $\Aut_\infty^{spin}(Q_M) \cong \Aut(\infty (-E_8) \oplus
\infty H)$ and one clearly sees that it is
independent of the initial manifold $M$.  This is not quite true for
$\Gamma_\infty^{spin}(M)$, but by Corollary \ref{twists-are-trivial2}
it does not depend on $M$ after localization away from 2.

By the proofs of Lemmas \ref{telescope} and \ref{telescope-equiv}, 
\begin{lemma}
  There is a homotopy equivalence $\mathcal{Y}_\infty^{spin}(n) \simeq
  \Z^2\times B\Aut(\infty (-E_8) \oplus \infty H)$, and the 2--functor
  $\calC_{spin} \to \calK_{spin}$
  induces a map $\mathcal{X}_\infty^{spin}(n) \to \mathcal{Y}_\infty^{spin}(n)$
  which is a homotopy equivalence away from the prime 2.
\end{lemma}

The group completion argument in the proof of Theorem
\ref{cat-equiv-spin} now goes through exactly as for the non-spin
version.

\section{Another infinite loop space operad}\label{operad-section}

In \cite{Tillmann-operad} Tillmann gives a general construction which
takes as input a family of group(oid) normal extensions of symmetric groups
equipped with appropriate wreath products.  The output of the
construction is an operad.  When the extensions are (stably)
homologically trivial Tillmann shows that spaces with an
action of the resulting operad are infinte loop spaces.

At the time of publication of \cite{Tillmann-operad} mapping class
groups of surfaces (and their variants) were the only known examples
of families of extensions which are \emph{not} already trivial at the
group level and which produce an infinte loop operad.  We observe now
that mapping class groups of simply connected $4$--manifolds also form
such a family.

Consider the family of groupoid extensions
\begin{equation}\label{groupoid-ext}
H_n \hookrightarrow G_n \twoheadrightarrow \Sigma_n
\end{equation}
where $H_n=\mathcal{C}(n,1)$ (see Definition \ref{CDef}), and $G_n$ is
the larger groupoid of isotopy classes of diffeomorphisms which
preserve the parametrization of the boundary but are no longer
required to preserve the ordering.  The epimorphism to the symemtric
group is given by sending an isotopy class to the permutation of the
in-coming boundary components that it induces.  Note that $\Sigma_n$
acts on $H_n$ by permuting the ordering of the in-coming boundary
components.

There are associative wreath products
\[
G_n \wr G_k \to G_{nk}
\]
induced by gluing the out-going boundary component of each of $n$
manifolds of type $(k,1)$ to the in-coming boundary components of a
manifold of type $(n,1)$.  The operad $\mathcal{E}$ formed from the
above family of extensions (\ref{groupoid-ext}) is
\[
\mathcal{E}_n := BG_n,
\]
and the operad composition maps are induced by the wreath products.
The component of $\mathcal{E}_1$ corresponding to the identity
morphism in $\mathcal{C}(1,1)$ is a point; this gives a unit for the
operad.  There is also a product given by a 4--sphere with three discs
removed, but this product is not strictly associative or unital.  One
could correct this with a quotient construction, but that is not
necessary for us.

The extensions (\ref{groupoid-ext}) are nontrivial.  However, they
become homologically trivial when stabilized (as with
$\mathcal{X}_\infty(n)$ in the previous section) by gluing copies of
$\C P^2 \# \overline{\C P^2} - \{\mbox{2 discs}\}$ to the in-coming
boundary component and extending isotopy classes by the identity.  The
resulting stablized extensions
\[
H_{\infty,n} \to G_{\infty,n} \to \Sigma_n
\]
have $H_{\infty,n} \cong \pi_1 \mathcal{X}_\infty(n) \cong
O_{\infty,\infty}(\Z)$, and $G_{\infty,n}$ the obvious analogue where
boundary components can be permuted.  Since $H_{\infty,n} \cong
H_{\infty,0}$ the action of $\Sigma_n$ on $H_{\infty,n}$, and hence on
$H_*(H_{\infty,n})$, is trivial.

There is also a spin analogue of the above discussion, leading to an
operad $\mathcal{E}^{spin}$.  Tillmann's argument
\cite{Tillmann-operad} applies verbatim to these.
\begin{theorem}
Algebras over the operad $\mathcal{E}$ (or $\mathcal{E}^{spin}$) are
canonically infinite loop spaces.
\end{theorem}

One should be able to adapt Wahl's comparison \cite{Wahl} to these
operads to show that the infinite loop space structure detected on the
stable 4--manifold (spin) mapping class group agrees with the usual
infinite loop space structure on $BO_{\infty,\infty}(\Z)^+$
(resp. $B\Aut(\infty (-E_8)\oplus \infty H)$).

\bibliographystyle{amsalpha} 
\bibliography{bib4man}

\end{document}